\documentclass[11pt]{article}

\title{Poisson , Dobinski , Rota and coherent states-
a fortieth   anniversary  memoir }
\author{A.K.Kwa\'sniewski\\  
\\ High School of Mathematics and Applied Informatics\\
PL - 15-021 Bialystok , ul.Kamienna 17,  Poland
\\e-mail: kwandr@uwb.edu.pl}

\usepackage{amsmath,amsthm}

\chardef\bslash=`\\ 
\begin{document}
\maketitle \ ArXiv: math. CO/0402125 v1 9 Feb 2004 corrected 17
Feb 2004

Forty years ago Rota G. C. \cite{1} while  proving that the
exponential generating function for Bell numbers $B_n$ is of the
form
\begin{equation}
                       \sum_{n=0}^\infty \frac {x^n}{n!}(B_n)= \exp(e^x-1)
\end{equation}
 used the linear functional  \textit{L } such that

\begin{equation}
                         L(X^{\underline{n}})=1  , \qquad  n\geq 0
\end{equation}
Then Bell numbers (see: formula (4)  in \cite{1})  are defined by
$$
            L(X^n)=B_n  ,\qquad   n \geq 0.
$$
Let us notice then that the above formula is exactly the Dobinski
formula \cite{2}  if   $L$  is interpreted as the average
functional for the random variable  $X$  with the Poisson
distribution with $L(X) = 1$.   (It is Blissard calculus inspired
umbral formula \cite{1}).\\
    Quite recently an interest to Stirling numbers and consequently
to Bell numbers was revived among "coherent states physicists"
\cite{3,4}. Namely the expectation value with respect to coherent
state $|\gamma >$ with $|\gamma| = 1$ of the $n$-th power of the
number of quanta operator is "just" the $n$-th Bell number $B_n$
and the explicit formula for this expectation number of quanta is
"just" Dobinski formula
\cite{3}.\\
 One faces the same situation with the $q$-coherent
states case \cite{3} i.e. the expectation value with respect to
$q$-coherent state $|\gamma> $ with  $|\gamma| = 1$ of the $n$-th
power of the number operator  is the $n$-th  $q$-Bell number
\cite{5} and the explicit formula becomes $q$-Dobinski formula.
Let us notice then that similar new  $q$-Dobinski  formula valid
for a new q-analogue of Stirling numbers of Cigler \cite{6} might
also be interpreted as the average of random variable $X_q^{n}$
with the same Poisson distribution with $L(X) = 1$ i.e.
\begin{equation}
                 L(X_q^n)=B_n(q)  ,\qquad n\geq 0 ; \qquad X_q^n\equiv X(X-1+q)...(X-1+q^{n-1}).
\end{equation}
For that to see use the identity by Cigler \cite{6}
\begin{equation}
 (x+1)(x+q)...(x+q^{n-1})=\sum_{k=0}^{n}\Big\{ {n \atop k}\Big\}_q  (x+1)^{\underline k}
\end{equation}
These Cigler $q$-analogue of Stirling numbers $\big\{ {n \atop
k}\big\}_q $ were given in \cite{6}  a combinatorial
interpretation in terms of weighted partitions. Therefore new
$q$-Bell numbers introduced above seem to deserve attention. In
\cite{7} a family of the so called $\psi$-Poisson processes was
introduced. The corresponding choice of the function sequence
$\psi$  leads to the $q$-Poisson process.

\end{document}